\documentclass[a4paper,11pt]{article}
\usepackage{amsmath,amsfonts,amsthm,dsfont,amssymb}
\usepackage[blocks]{authblk}
\usepackage{cite}
\usepackage{graphicx}
\usepackage{amssymb,array,arydshln}
\usepackage{amsfonts}
\usepackage{amsbsy}
\usepackage{cite}
\usepackage{amssymb}
\usepackage{amsmath}
\usepackage{graphics}
\newenvironment{keywords}{\noindent\textbf{Keywords:}}{}
\newenvironment{classification}{\noindent\textbf{AMS subject classifications.}}{}
\date{}
\newcommand{\email}[1]{\texttt{\small #1}}
\newtheorem{theorem}{Theorem}[section]
\newtheorem{remark}[theorem]{Remark}
\newtheorem{example}[theorem]{Example}
\newtheorem{lemma}[theorem]{Lemma}
\newtheorem{corollary}[theorem]{Corollary}
\newtheorem{definition}[theorem]{Definition}
\newtheorem{proposition}[theorem]{Proposition}

\newtheorem{conjecture}{Conjecture}
\newtheorem{question}{Question}
\newcommand{\bea}{\begin{eqnarray}}
\newcommand{\eea}{\end{eqnarray}}
\newcommand{\beq}{\begin{eqnarray*}}
\newcommand{\eeq}{\end{eqnarray*}}
\def \bd{\begin{definition}}
\def \ed{\end{definition}}
\def \bqu{\begin{question}}
\def \equ{\end{question}}
\def \bcc{\begin{conjecture}}
\def \ecc{\end{conjecture}}
\def \bt{\begin{theorem}}
\def \et{\end{theorem}}
\def \bl{\begin{lemma}}
\def \el{\end{lemma}}
\def \bc{\begin{corollary}}
\def \ec{\end{corollary}}
\def \be{\begin{equation}}
\def \ee{\end{equation}}
\def \ben{\begin{enumerate}}
\def \een{\end{enumerate}}
\def \ba{\begin{array}}
\def \ea{\end{array}}
\def \bp{\begin{proposition}}
\def \ep{\end{proposition}}
\def \bx{\begin{example}}
\def \ex{\end{example}}
\def \br{\begin{remark}}
\def \er{\end{remark}}
\def \bdsc{\begin{description}}
\def \edsc{\end{description}}

\def\1{1\!\!1}
\def\0{0\!\!0}

\begin{document}
	\title{On adjacency and (signless) Laplacian spectra of centralizer and co-centralizer graphs of some finite non-abelian groups
	}

\author[1]{Jharna Kalita}
\author[2]{Somnath Paul\footnote{Corresponding  Author.}}
\affil[1]{Department of Applied Sciences\\ Tezpur University\\ Napaam-784028, Assam, India. \email{app21104@tezu.ac.in}}
\affil[2]{Department of Applied Sciences\\ Tezpur University\\ Napaam-784028, Assam, India. \email{som@tezu.ernet.in}}

\pagestyle{myheadings} \markboth{J. Kalita \& S. Paul}{On adjacency and (signless)laplacian spectra of centralizer $\ldots$}
   \maketitle
\begin{abstract}
Let $G$ be a finite non abelian group. The centralizer graph of $G$  is a simple undirected graph $\Gamma_{cent}(G)$, whose vertices are the proper centralizers of $G$ and two vertices are adjacent if and only if their cardinalities are identical {\rm\cite{omer}}. The complement of the centralizer graph is called the co-centralizer graph. In this paper, we investigate the adjacency and (signless) Laplacian spectra of centralizer and co-centralizer graphs of some classes of finite non-abelian groups and obtain some conditions on a group so that the centralizer and co-centralizer graphs are adjacency, (signless) Laplacian integral.
\end{abstract}


\section{Introduction}
Let $G$ be a finite non-ableian group. In literature, there are many occasions when one associates a graph to a group $G$ in different ways. For example, a commuting graph is a graph with $G\setminus Z(G)$ as the vertex set and two vertices $x$ and $y$ are adjacent if and only if $xy=yx.$ Similarly, the non-commuting graph is a graph with $G\setminus Z(G)$ as the set of vertices and two vertices $x$ and $y$ are adjacent if and only if $xy\ne yx.$ Likewise, non-nilpotent graph, cyclic, non-cyclic and conjugacy class graphs has also been associated on a group. The centralizer graph of $G$  is a simple undirected graph $\Gamma_{cent}(G)$, whose vertices are the proper centralizers of $G$ and two vertices are adjacent if their cardinalities are identical {\rm\cite{omer}}. A brief study is done about the structure of centralizer graph in {\rm\cite{omer}}. Another definition is given for centralizer graph in {\rm\cite{gazvin}} where the centralizer graph of $G$ is a simple undirected graph with the proper centralizers of $G$ constituting the vertex set and two vertices are adjacent if they are same, and its complement graph is called non-centralizer graph. In this article, we consider the centralizer graph defined in {\rm\cite{omer}}. Also, we define the complement of that graph as the co-centralizer graph and denote it by $\overline{\Gamma_{cent}(G)}$ .\\
\indent \hspace{.5cm} For a simple graph $H$ on $n$ vertices, the adjacency matrix $A(H)$ is a matrix of order $n,$ whose  $(i,j)$-th entry is 1, if the $i$-th vertex is adjacent to the $j$-th vertex; otherwise it is 0. Also, the Laplacian (resp. signless Laplacian) matrix of $H$ is defined as $L(H)=D(H)-A(H)$ (resp. $Q(H)=D(H)+A(H)$), where $D(H)$ is the diagonal matrix of order $n,$ with degree of the $i$-th vertex as the $i$-th diagonal entry.

\indent \hspace{.5cm}If $M$ is a symmetric matrix, then the characteristic polynomial of $M$ has only real zeroes. We will represent this family of eigenvalues (known as the \textit{spectrum}) as $$\sigma_M=\left(
                  \begin{array}{cccc}
                    \mu_1 & \mu_2 &\cdots&\mu_p\\
                    m_1 & m_2 &\cdots & m_p
                  \end{array}
                \right),
$$ where $\mu_1,\mu_2,\ldots,\mu_p$ are the distinct eigenvalues of $M$ and $m_1, m_2,\ldots,m_p$ are the corresponding multiplicities. Since each of $A(H),~L(H)$ and $Q(H)$ is symmetric, we will refer the corresponding spectrum as  \textit{the adjacency,  the Laplacian} and \textit{the signless Laplacian spectrum}, respectively.

In \cite{rkn}, the adjacency spectrum of the commutating graph of some finite non-abelian groups is discussed. In \cite{jr}, the Laplacian and signless Laplacian spectra of the commutating graph of some finite non-abelian groups is investigated, whereas the Laplacian spectrum of the non-commutating graph of some finite non-abelian groups is determined in \cite{par}. In \cite{pir}, the Laplacian spectrum of unitary cayley graphs are discussed. For other related results the reader can look into \cite{rkn,jr,par,pir} and the references therein.

A graph is called adjacency (respectively (signless) Laplacian) integral if the adjacency (respectively (signless) Laplacian) spectrum consists entirely of integers. In this article, we consider some finite non-abelian groups, namely the generalized quaternion group, the dihedral group, the quasidihedral group, the metacyclic group, and the projective special linear group,  and investigate the adjacency, (signless) Laplacian spectra of centralizer and co-centralizer graphs of them. Moreover, we obtain some conditions so that their centralizer and co-centralizer graphs are adjacency, (signless) Laplacian integral.

\section{Preliminaries}\label{men}

Consider $Q_{4n} =<x,y : x^{2n} = 1, x^n = y^2, yx=x^{-1}y>,$ be the generalized quaternion group of order $4n,$  where $n\ge2$ , and $Z(Q_{4n})$= $\{1,x^n\}$. Then $Q_{4n}$ can be written as $A\cup B,$ where $A=\{1,x,x^2,\ldots,x^{2n-1}\}$ and $B=\{y,xy,x^2y,\ldots,x^{2n-1}y\},$ where each element of $B$ is of order 4. It has $n+1$ distinct centralizers with one of cardinality $2n$ and others are of cardinality $4$. For illustration,  we note that for any $z\in Z(Q_{4n}),$ and $1\le i\le 2n-1,$
\begin{eqnarray*}
  C_{Q_{4n}}(x)
  &=&C_{Q_{4n}}(x^iz)\\ &=& Z(Q_{4n})\cup x Z(Q_{4n})\cup x^2Z(Q_{4n})\cup\ldots\cup x^{n-1}Z(Q_{4n})\\
    &=& \{1,x^n\} \cup x \{1,x^n\}\cup x^2\{1,x^n\}\cup\ldots\cup x^{n-1}\{1,x^n\}\\
    &=& \{1,x^n\} \cup \{x,x^{n+1}\}\cup \{x^2,x^{n+2}\}\cup\ldots\cup \{x^{n-1},x^{2n-1}\} \\
    &=& \{1,x,x^2,\ldots,x^{2n-1}\}.
\end{eqnarray*}
Moreover, for $1\le j\le n,$
\begin{eqnarray*}
  C_{Q_{4n}}(yx^j)=C_{Q_{4n}}(yx^jz) &=& Z(Q_{4n})\cup yx^j Z(Q_{4n}) \\
    &=& \{1,x^n\} \cup yx^j \{1,x^n\}\\
    &=& \{1,x^n\} \cup \{yx^j,yx^{n+j}\}.
\end{eqnarray*}

Therefore, from the definition of centralizer graph, it follows that $\Gamma_{cent}(Q_{4n})$ is a graph with $n+1$ vertices where one component is $K_n$ and the other is an isolated vertex, i.e., $\Gamma_{cent}(Q_{4n})\cong K_n\sqcup K_1$. As co-centralizer graph is the complement of this graph, therefore $\overline{\Gamma_{cent}(Q_{4n})}\cong K_{1,n}$.

The following result gives the adjacency characteristics polynomial for $K_{p_1,p_2,\ldots,p_n}$ and will be useful to derive some of our main results.
\bl{\rm \cite{spec}}\label{lem1} The adjacency characteristics polynomial of the complete multipartite graph $K_{p_1,p_2,\ldots,p_n}$ , where $p_1+p_2+ \ldots+p_n=P$ is:

\begin{equation}\label{e1}
 P_{A(G)}(\lambda)=\lambda^{P-n}\left[\prod_{i=1}^n(\lambda+p_i)-\sum_{i=1}^n p_i\prod_{j=1,j\neq i}^n(\lambda+p_j)\right].
\end{equation} \el

\section{Spectra of centralizer graphs of some finite non-abelian groups}
\subsection{Spectra of $\Gamma_{cent}(Q_{4n})$}\label{sec1}
In this section, we consider the centralizer graph of $Q_{4n}$ and obtain the adjacency, Laplacian and signless Laplacian spectra of it. It is well known (see \cite{spg}) that if  $G=K_{m_1}\sqcup K_{m_2}\sqcup \ldots\sqcup K_{m_l}$, then  $$\sigma(A(G))=\left(
                       \begin{array}{ccccc}
                         -1 & m_1-1 & m_2-1 & \ldots & m_l-1 \\
                         \displaystyle\sum_{i=1}^l{m_i}-l & 1 & 1 & \ldots & 1 \\
                       \end{array}
                     \right).$$

Similarly, if $G=l_1K_{m_1}\sqcup l_2K_{m_2}\sqcup \ldots\sqcup l_kK_{m_k},$ then $$\sigma(L(G))=\left(
                       \begin{array}{ccccc}
                          0 & m_1 & m_2 & \ldots & m_k \\
                         \displaystyle\sum_{i=1}^k{l_i} & l_1(m_1-1) & l_2(m_2-1)  & \ldots & l_k(m_k-1) \\
                       \end{array}
                     \right).$$

Therefore, $\sigma(A(\Gamma_{cent}(Q_{4n}))=\left(
                       \begin{array}{ccccc}
                         -1 & 0 & n-1 \\
                         n-1 & 1 & 1  \\
                       \end{array}
                     \right)$
and $\sigma(L(\Gamma_{cent}(Q_{4n}))=\left(
                       \begin{array}{ccccc}
                         0 &  n \\
                         2 &  n-1 \\
                       \end{array}
                     \right)$. It is easily observed that  $\sigma(Q(\Gamma_{cent}(Q_{4n}))=\left(
                       \begin{array}{ccccc}
                         0 & n-2 & 2(n-1) \\
                         1 & n-1 & 1  \\
                       \end{array}
                     \right)$.\\
Thus, $\Gamma_{cent}(Q_{4n})$ is adjacency, Laplacian and signless Laplacian integral for any $n$.
\subsection{Spectra of $\Gamma_{cent}(D_{2n})$}
In this section, we consider the centralizer graph of the dihedral group $D_{2n}=<x,y : x^n = y^2 = 1, yxy^{-1}=x^{-1}>$ and obtain the adjacency, Laplacian and signless Laplacian spectra of it. The centralizer graph $\Gamma_{cent}(D_{2n})$ of dihedral group $D_{2n}$ is $K_1\sqcup K_n,$ when $n$ is odd, and is $K_1\sqcup K_{\frac{n}{2}},$ when $n$ is even.

Therefore,  if $n$ is odd, then $\sigma(A(\Gamma_{cent}(D_{2n}))=\left(
                       \begin{array}{ccccc}
                         -1 & 0 & n-1 \\
                         n-1 & 1 & 1  \\
                       \end{array}
                     \right),$ $\sigma(L(\Gamma_{cent}(D_{2n}))=\left(
                       \begin{array}{ccccc}
                         0 & n \\
                         2 & n-1  \\
                       \end{array}
                     \right)$ and $\sigma(Q(\Gamma_{cent}(D_{2n}))=\left(
                       \begin{array}{ccccc}
                         0 & n-2 & 2(n-1) \\
                         1 & n-1 & 1  \\
                       \end{array}
                     \right).$ Also, if $n$ is even, then $\sigma(A(\Gamma_{cent}(D_{2n}))=\left(
                       \begin{array}{ccccc}
                         -1 & 0 & \frac{n}{2}-1 \\
                         \frac{n}{2}-1 & 1 & 1  \\
                       \end{array}
                     \right),$  \linebreak$\sigma(L(\Gamma_{cent}(D_{2n}))=\left(
                       \begin{array}{ccccc}
                          0 & \frac{n}{2} \\
                          2 & \frac{n}{2}-1  \\
                       \end{array}
                     \right),$ and $\sigma(Q(\Gamma_{cent}(D_{2n}))=\left(
                       \begin{array}{ccccc}
                         0 & \frac{n}{2}-2 & n-2 \\
                         1 & \frac{n}{2}-1 & 1  \\
                       \end{array}
                     \right).$

Thus, $\Gamma_{cent}(D_{2n})$ is adjacency, Laplacian and signless Laplacian integral for any $n$.\\
\br\label{rm1} Let us consider the metacyclic group $M_{2pq}=<a,b: a^p=b^{2q}=1,bab^{-1}=a^{-1}>$ , where $p>2.$ It can be easily observed that for even $p$ (respectively for odd $p$) the corresponding centralizer graph is same as that of the centralizer graph of dihedral group $D_{2n}$ for even $n$ (respectively for odd $n$), and  is independent of $q.$ Therefore the adjacency, Laplacian and signless Laplacian spectra of $M_{2pq}$ is exactly same as the corresponding spectra of $D_{2n}$.\er

\subsection{Spectra of $\Gamma_{cent}(QD_{2^n})$}
In this section, we consider the centralizer graph of the quasidihedral group $QD_{2^n}=<a,b:a^{2^{n-1}}=b^2=1,bab^{-1}=a^{2^{n-2}-1}>,$ where $n\ge4,$ and obtain the adjacency, Laplacian and signless Laplacian spectra of it. The centralizer graph $\Gamma_{cent}(QD_{2^n})$ of $QD_{2^n}$ is $K_1\sqcup K_{2^{n-2}}$.

Therefore, $\sigma(A(\Gamma_{cent}(Q_{4n}))=\left(
                       \begin{array}{ccccc}
                         -1 & 0 & 2^{n-2}-1 \\
                         2^{n-2}-1 & 1 & 1  \\
                       \end{array}
                     \right),$\linebreak $\sigma(L(\Gamma_{cent}(QD_{2^n}))=\left(
                       \begin{array}{ccccc}
                         0 &  2^{n-2} \\
                         2 &  2^{n-2}-1 \\
                       \end{array}
                     \right),$ and \linebreak $\sigma(Q(\Gamma_{cent}(QD_{2^n}))=\left(
                       \begin{array}{ccccc}
                         0 & 2^{n-2}-2 & 2^{n-1}-2 \\
                         1 & 2^{n-2}-1 & 1  \\
                       \end{array}
                     \right)$.

Thus, $\Gamma_{cent}(QD_{2^n})$ is adjacency, Laplacian and signless Laplacian integral for any $n$.
\subsection{Spectra of $\Gamma_{cent}(PSL(2,2^k))$}\label{m2}
In this section, we consider the centralizer graph of the projective special linear group $PSL(2,2^k)$ and obtain the adjacency, Laplacian and signless Laplacian spectra of it. The centralizer graph $\Gamma_{cent}(PSL(2,2^k))$ of $PSL(2,2^k)$ is $K_{2^k+1}\sqcup K_{2^{k-1}(2^k+1)}\sqcup K_{2^{k-1}(2^k-1)}$.

Therefore, $$\sigma(A(\Gamma_{cent}(PSL(2,2^k)))=\left(
                       \begin{array}{ccccc}
                         -1 & 2^k &  2^{k-1}(2^k+1)-1  &  2^{k-1}(2^k-1)-1  \\
                         2^{2k}+2^{k}-2 & 1 & 1 & 1  \\
                       \end{array}
                     \right),$$ and $\sigma(L(\Gamma_{cent}(PSL(2,2^k)))=\left(
                       \begin{array}{ccccc}
                         0 &  2^k+1 & 2^{k-1}(2^k+1) & 2^{k-1}(2^k-1)\\
                         3 &  2^k &  2^{k-1}(2^k+1)-1  &  2^{k-1}(2^k-1)-1 \\
                       \end{array}
                     \right)$.

Let $\1_n$ (resp. $\0_n$) denote the $n\times 1$ vector with each entry 1 (resp. 0). Also, let $J_n$ (resp. $\textbf{0}_n$) denote the matrix of order $n$ with all entries equal to 1 (resp. 0) (we will write $J$ (resp. \textbf{0})  if the order is clear from the context).   The following theorem describes the signless Laplacian spectrum of $\Gamma_{cent}(PSL(2,2^k)).$

\bt\label{dqt1} Let $\Gamma_{cent}(PSL(2,2^k))$ be the centralizer graph of the projective special linear group. Then
\begin{enumerate}
\item[{\rm(a)}] $2^k-1\in\sigma(Q(\Gamma_{cent}(PSL(2,2^k))))$ with multiplicity $2^k;$
\item[{\rm(b)}] $2^{k-1}(2^k+1)-2\in\sigma( Q(\Gamma_{cent}(PSL(2,2^k))))$ with multiplicity $2^{k-1}(2^k+1)-1;$
\item[{\rm(c)}] $2^{k-1}(2^k-1)-2\in\sigma( Q(\Gamma_{cent}(PSL(2,2^k))))$  with multiplicity $2^{k-1}(2^k-1)-1,$
\item[{\rm(d)}] $(2^{k}+1)(2^k-2)\in\sigma( Q(\Gamma_{cent}(PSL(2,2^k))))$  with multiplicity $1,$
\item[{\rm(e)}] $2^{2k}+2^k-2\in\sigma( Q(\Gamma_{cent}(PSL(2,2^k))))$  with multiplicity $1,$
\item[{\rm(f)}] $2^{k+1}\in\sigma( Q(\Gamma_{cent}(PSL(2,2^k))))$  with multiplicity $1.$
\end{enumerate}\et
{\bf Proof.}  With a suitable labeling of the vertices, the signless Laplacian matrix for $\Gamma_{cent}(PSL(2,2^k))$ can be written as
\begin{eqnarray*}
  &&Q(\Gamma_{cent}(PSL(2,2^k))) \\
  &=& \left(
                                   \begin{array}{c|c|c}
                                     J+(2^k-1)I & \textbf{0} & \textbf{0} \\
                                     \hline
                                     \textbf{0} & J+(2^{k-1}(2^k+1)-2)I & \textbf{0} \\
                                     \hline
                                     \textbf{0} & \textbf{0} & J+(2^{k-1}(2^k-1)-2)I\\
                                   \end{array}
                                 \right).
\end{eqnarray*}

Now, $Q(\Gamma_{cent}(PSL(2,2^k)))\left(
           \begin{array}{c}
            -1\\
             1\\
             \0_{2^k-1}\\
             \hline
             \0_{2^{k-1}(2^k+1)}\\
             \hline
             \0_{2^{k-1}(2^k-1)}\\
           \end{array}
         \right)=(2^k-1)\left(
           \begin{array}{c}
             -1\\
             1\\
             \0_{2^k-1}\\
             \hline
             \0_{2^{k-1}(2^k+1)}\\
             \hline
             \0_{2^{k-1}(2^k-1)}\\
           \end{array}
         \right).$

Therefore, $(2^k-1)$ is an eigenvalue of $Q(\Gamma_{cent}(PSL(2,2^k))),$ and the following set $V_1$ lists the set of $2^k$ linearly independent eigenvectors corresponding to the eigenvalue $2^k-1;$ $$V_1=\left\{ \left(
           \begin{array}{c}
            -1\\
             1\\
             \0_{2^k-1}\\
             \hline
             \0_{2^{k-1}(2^k+1)}\\
             \hline
             \0_{2^{k-1}(2^k-1)}\\
           \end{array}
         \right),  \left(
           \begin{array}{c}
            -1\\
            0\\
             1\\
             \0_{2^k-2}\\
             \hline
             \0_{2^{k-1}(2^k+1)}\\
             \hline
             \0_{2^{k-1}(2^k-1)}\\
           \end{array}
         \right), \ldots, \left(
           \begin{array}{c}
            -1\\
             \0_{2^k-1}\\
             1\\
             \hline
             \0_{2^{k-1}(2^k+1)}\\
             \hline
             \0_{2^{k-1}(2^k-1)}\\
           \end{array}
         \right)\right\}.$$
Again, $Q(\Gamma_{cent}(PSL(2,2^k)))\left(
           \begin{array}{c}
             \0_{2^k+1}\\
             \hline
             -1 \\
             1 \\
             \0_{2^{k-1}(2^k+1)-2}\\
             \hline
             \0_{2^{k-1}(2^k-1)}\\
           \end{array}
         \right)=(2^{k-1}(2^k+1)-2)\left(
           \begin{array}{c}
             \0_{2^k+1}\\
             \hline
             -1 \\
             1 \\
             \0_{2^{k-1}(2^k+1)-2}\\
             \hline
             \0_{2^{k-1}(2^k-1)}\\
           \end{array}
         \right).$
Therefore, $(2^{k-1}(2^k+1)-2)$ is an eigenvalue of $Q(\Gamma_{cent}(PSL(2,2^k))),$ and the set $V_2$ gives $2^k$ linearly independent eigenvectors corresponding to it; $$V_2=\left\{\left(
           \begin{array}{c}
             \0_{2^k+1}\\
             \hline
             -1 \\
             1 \\
             \0_{2^{k-1}(2^k+1)-2}\\
             \hline
             \0_{2^{k-1}(2^k-1)}\\
           \end{array}
         \right),\left(
           \begin{array}{c}
             \0_{2^k+1}\\
             \hline
             -1 \\
             0\\
             1 \\
             \0_{2^{k-1}(2^k+1)-3}\\
             \hline
             \0_{2^{k-1}(2^k-1)}\\
           \end{array}
         \right), \ldots, \left(
           \begin{array}{c}
             \0_{2^k+1}\\
             \hline
             -1 \\
             \0_{2^{k-1}(2^k+1)-2}\\
             1 \\
             \hline
             \0_{2^{k-1}(2^k-1)}\\
           \end{array}
         \right)\right\}.$$

Similarly, $$Q(\Gamma_{cent}(PSL(2,2^k)))\left(
           \begin{array}{c}
             \0_{2^k+1}\\
             \hline
             \0_{2^{k-1}(2^k+1)}\\
             \hline
             -1 \\
             1 \\
             \0_{2^{k-1}(2^k-1)-2}\\
           \end{array}
         \right)=(2^{k-1}(2^k-1)-2)\left(
           \begin{array}{c}
             \0_{2^k+1}\\
             \hline
             \0_{2^{k-1}(2^k+1)}\\
             \hline
             -1 \\
             1 \\
             \0_{2^{k-1}(2^k-1)-2}\\
           \end{array}
         \right).$$

Therefore, $(2^{k-1}(2^k-1)-2)$ is an eigenvalue of $Q(\Gamma_{cent}(PSL(2,2^k))),$ and the set set $V_3$ gives $(2^{k-1}(2^k-1)-1)$ linearly independent eigenvectors corresponding to it;
$$V_3=\left\{\left(
           \begin{array}{c}
             \0_{2^k+1}\\
             \hline
             \0_{2^{k-1}(2^k+1)}\\
             \hline
             -1 \\
             1 \\
             \0_{2^{k-1}(2^k-1)-2}\\
           \end{array}
         \right),\left(
           \begin{array}{c}
             \0_{2^k+1}\\
             \hline
             \0_{2^{k-1}(2^k+1)}\\
             \hline
             -1 \\
             0 \\
             1 \\
             \0_{2^{k-1}(2^k-1)-3}\\
           \end{array}
         \right), \ldots, \left(
           \begin{array}{c}
             \0_{2^k+1}\\
             \hline
             \0_{2^{k-1}(2^k+1)}\\
             \hline
             -1 \\
             \0_{2^{k-1}(2^k-1)-2}\\
             1 \\
           \end{array}
         \right)\right\}.$$

 Moreover,
 \begin{eqnarray*}
   Q(\Gamma_{cent}(PSL(2,2^k)))\left(
           \begin{array}{c}
             \0_{2^k+1}\\
             \hline
             \0_{2^{k-1}(2^k+1)}\\
             \hline
             \1_{2^{k-1}(2^k-1)}\\
           \end{array}
         \right)&=&(2^k+1)(2^k-2)\left(
           \begin{array}{c}
             \0_{2^k+1}\\
             \hline
             \0_{2^{k-1}(2^k+1)}\\
             \hline
             \1_{2^{k-1}(2^k-1)}\\
           \end{array}
         \right), \\
   Q(\Gamma_{cent}(PSL(2,2^k)))\left(
           \begin{array}{c}
             \0_{2^k+1}\\
             \hline
             \1_{2^{k-1}(2^k+1)}\\
             \hline
             \0_{2^{k-1}(2^k-1)}\\
           \end{array}
         \right)&=&(2^{2k}+2^k-2)\left(
           \begin{array}{c}
             \0_{2^k+1}\\
             \hline
             \1_{2^{k-1}(2^k+1)}\\
             \hline
             \0_{2^{k-1}(2^k-1)}\\
           \end{array}
         \right), \\
   Q(\Gamma_{cent}(PSL(2,2^k)))\left(
           \begin{array}{c}
             \1_{2^k+1}\\
             \hline
             \0_{2^{k-1}(2^k+1)}\\
             \hline
             \0_{2^{k-1}(2^k-1)}\\
           \end{array}
         \right)&=&(2^{k+1})\left(
           \begin{array}{c}
             \1_{2^k+1}\\
             \hline
             \0_{2^{k-1}(2^k+1)}\\
             \hline
             \0_{2^{k-1}(2^k-1)}\\
           \end{array}
         \right).
 \end{eqnarray*}

 Therefore, $(2^k+1)(2^k-2),(2^{2k}+2^k-2),$ and $(2^{k+1})$ are an eigenvalues of $Q(\Gamma_{cent}(PSL(2,2^k)))$ with $\left(
           \begin{array}{c}
             \0_{2^k+1}\\
             \hline
             \0_{2^{k-1}(2^k+1)}\\
             \hline
             \1_{2^{k-1}(2^k-1)}\\
           \end{array}
         \right), \left(
           \begin{array}{c}
             \0_{2^k+1}\\
             \hline
             \1_{2^{k-1}(2^k+1)}\\
             \hline
             \0_{2^{k-1}(2^k-1)}\\
           \end{array}
         \right),$ and $\left(
           \begin{array}{c}
             \1_{2^k+1}\\
             \hline
             \0_{2^{k-1}(2^k+1)}\\
             \hline
             \0_{2^{k-1}(2^k-1)}\\
           \end{array}
         \right)$ as the corresponding eigenvector, respectively.

We note that $V_1\cup V_2 \cup V_3 \cup \left(
           \begin{array}{c}
             \0_{2^k+1}\\
             \hline
             \0_{2^{k-1}(2^k+1)}\\
             \hline
             \1_{2^{k-1}(2^k-1)}\\
           \end{array}
         \right) \cup \left(
           \begin{array}{c}
             \0_{2^k+1}\\
             \hline
             \1_{2^{k-1}(2^k+1)}\\
             \hline
             \0_{2^{k-1}(2^k-1)}\\
           \end{array}
         \right) \cup \left(
           \begin{array}{c}
             \1_{2^k+1}\\
             \hline
             \0_{2^{k-1}(2^k+1)}\\
             \hline
             \0_{2^{k-1}(2^k-1)}\\
           \end{array}
         \right)$ is a set of mutually orthogonal eigenvectors for $\Gamma_{cent}(PSL(2,2^k)).$ Since the order of $\Gamma_{cent}(PSL(2,2^k))$ is $2^{2k}+2^k+1$, the result follows.\qed

Thus, $\Gamma_{cent}(PSL(2,2^k))$ is adjacency, Laplacian and signless Laplacian integral for any $k$.
\section{Spectra of co-centralizer graphs of some finite non-abelian groups}

\subsection{Spectra of $\overline{\Gamma_{cent}(Q_{4n})}$}\label{sec2}

In this section, we consider the co-centralizer graph of $Q_{4n}$ and obtain the adjacency, Laplacian and signless Laplacian spectra of it. It is well known (see \cite{spg}) that   the adjacency spectra of a complete bipartite graph $K_{m,n}$ is $\left(
                                                               \begin{array}{ccc}
                                                                 \sqrt{mn} & -\sqrt{mn} & 0 \\
                                                                 1 & 1 & m+n-2 \\
                                                               \end{array}
                                                             \right).$
 As it is already observed in Section~\ref{men},  $\overline{\Gamma_{cent}(Q_{4n})}=K_{1,n}$. Therefore, $\sigma(A(\overline{\Gamma_{cent}(Q_{4n}})))=\left(\begin{array}{ccc}
                                                                                                                                           \sqrt{n} & -\sqrt{n} & 0\\
                                                                                                                                           1 & 1 & n-1
                                                                                                                                         \end{array}\right).$
 Therefore,  $\overline{\Gamma_{cent}(Q_{4n})}$ is adjacency integral if $n$ is a perfect square. Also,  by Lemma~5 of {\cite{pir}}, $\sigma(L(\overline{\Gamma_{cent}(Q_{4n}})))=\left(
                                                                                            \begin{array}{ccc}
                                                                                              0 & 1 & 1+n \\
                                                                                              1 & n-1 & 1 \\
                                                                                            \end{array}
                                                                                          \right).$ Since, for a bipartite graph the Laplacian spectrum coincides with the signless Laplacian spectrum, we have $\sigma(L(\overline{\Gamma_{cent}(Q_{4n}})))=\sigma(Q(\overline{\Gamma_{cent}(Q_{4n}}))).$  Thus, $\overline{\Gamma_{cent}(Q_{4n})}$ is Laplacian and signless Laplacian integral for any value of $n.$
\subsection{Spectra of $\overline{\Gamma_{cent}(D_{2n})}$}
In this section, we consider the co-centralizer graph of dihedral group $D_{2n}$ and obtain the adjacency, Laplacian and signless Laplacian spectra of it. The co centralizer graph of $D_{2n}$ is $\overline{\Gamma_{cent}(D_{2n})}=\left\{
                                                                                           \begin{array}{ll}
                                                                                             K_{1,n}, & \hbox{if $n$ is odd;} \\
                                                                                             K_{1,\frac{n}{2}}, & \hbox{if $n$ is even.}
                                                                                           \end{array}
                                                                                         \right.$

Therefore, $\sigma(A(\overline{\Gamma_{cent}(D_{2n})}))=\left\{
                                                          \begin{array}{ll}
                                                            \left(
                                                                           \begin{array}{ccc}
                                                                             \sqrt{n} & -\sqrt{n} & 0 \\
                                                                             1 & 1 & n-1 \\
                                                                           \end{array}
                                                                         \right), & \hbox{if $n$ is odd;} \\
                                                            \left(
                                                                           \begin{array}{ccc}
                                                                             \sqrt{\frac{n}{2}} & -\sqrt{\frac{n}{2}} & 0 \\
                                                                             1 & 1 & \frac{n}{2}-1 \\
                                                                           \end{array}
                                                                         \right), & \hbox{if $n$ is even.}
                                                          \end{array}
                                                        \right.$
Thus, $\overline{\Gamma_{cent}(D_{2n})}$ is adjacency integral if $n$ is a perfect square for odd $n,$ and $\frac{n}{2}$ is a perfect square for even $n.$

As discussed in Subsection~\ref{sec2}, it can be seen that $$\sigma(L(\overline{\Gamma_{cent}(D_{2n})}))=\sigma(Q(\overline{\Gamma_{cent}(D_{2n})}))=\left\{
                                                          \begin{array}{ll}
                                                            \left(
                                                                           \begin{array}{ccc}
                                                                            0 & 1+n & 1 \\
                                                                             1 & 1 & n-1 \\
                                                                           \end{array}
                                                                         \right), & \hbox{if $n$ is odd;} \\
                                                            \left(
                                                                           \begin{array}{ccc}
                                                                             0 & 1+\frac{n}{2} & 1 \\
                                                                             1 & 1 & \frac{n}{2}-1 \\
                                                                           \end{array}
                                                                         \right), & \hbox{if $n$ is even.}
                                                          \end{array}
                                                        \right.$$

Hence for any value of $n,$ $\overline{\Gamma_{cent}(D_{2n})})$ is both Laplacian and signless Laplacian integral.

\br By virtue of Remark~\ref{rm1}, we can conclude that the adjacency, Laplacian and signless Laplacian spectrum of $\overline{\Gamma_{cent}(M_{2pq})}$ is exactly same as the corresponding spectrum of $\overline{\Gamma_{cent}(D_{2n})}$.\er

\subsection{Spectra of Quasidihedral group $\overline{\Gamma_{cent}(QD_{2^n})}$}

In this section, we consider the co-centralizer graph of the Quasidihedral group $QD_{2^n},$ where $n\ge4,$ and obtain the adjacency, Laplacian and signless Laplacian spectra of it. Since $\overline{\Gamma_{cent}(QD_{2^n})}=K_{1,2^{n-2}},$ it follows that $$\sigma(A(\overline{\Gamma_{cent}(QD_{2^n})}))=\left(
                                                                                         \begin{array}{ccc}
                                                                                          \sqrt{2^{n-2}} & -\sqrt{2^{n-2}} & 0 \\
                                                                                           1 & 1 & 2^{n-2}-1 \\
                                                                                         \end{array}
                                                                                       \right).$$
Thus, $\overline{\Gamma_{cent}(QD_{2^n})}$ is adjacency integral, if $2^{n-2}$ is a perfect square. Also, $\sigma(L(\overline{\Gamma_{cent}(QD_{2^n})}))=\sigma(Q(\overline{\Gamma_{cent}(QD_{2^n})}))=\left(
                                                                                         \begin{array}{ccc}
                                                                                          2^{n-2}+1 & 0 & 1 \\
                                                                                           1 & 1 & 2^{n-2}-1 \\
                                                                                         \end{array}
                                                                                       \right),$  showing that for any value of $n,$ $\overline{\Gamma_{cent}(QD_{2^n})})$ is both Laplacian and signless Laplacian integral.

\subsection{Spectra of $\overline{\Gamma_{cent}(PSL(2,2^k))}$}
As observed in Subsection~\ref{m2}, $\overline{\Gamma_{cent}(PSL(2,2^k))}$ is the complete tripartite graph $K_{2^k+1,2^{k-1}(2^k-1),2^{k-1}(2^k+1)}$. Therefore, by equation (\ref{e1})  we get,
\begin{eqnarray*}
P_{A(\overline{\Gamma_{cent}(PSL(2,2^k))}}(\lambda)&=&\lambda^{2^k+2^{2k}-2}[\lambda^3-\{2^{4k-2}+2^{3k}+3\times2^{2k-2}\}\lambda+\\
&&(-2^{5k-1}-2^{4k-1}+2^{3k-1}+2^{2k-1})].
\end{eqnarray*}

Hence, we have the following theorem which describes the adjacency spectrum of $\overline{\Gamma_{cent}(PSL(2,2^k))}.$

\bt\label{dqt1} Let $\overline{\Gamma_{cent}(PSL(2,2^k))}$ be the co-centralizer graph of the projective special linear group. Then $\sigma(A(\overline{\Gamma_{cent}(PSL(2,2^k))}))$ consists of
\begin{enumerate}
\item[{\rm(a)}] $0$ with multiplicity $2^k+2^{2k}-2;$
\item[{\rm(b)}] three roots of the equation $x^3-(2^{4k-2}+3(2^{2k-2})+2^{3k})x +(-2^{5k-1}-2^{4k-1}+2^{3k-1}+2^{2k-1})=0$.
\end{enumerate}\et

Also,
\begin{eqnarray*}
  &&\sigma(L(\overline{\Gamma_{cent}(PSL(2,2^k))})) \\
  &=& \left(
                                                                                                                                            \begin{array}{ccccc}
                                                                                                                                               0& 2^{2k} & 2^{k-1}+2^{2k-1}+1 & 2^{2k-1}+3(2^{k-1})+1 & 2^{2k}+2^k+1 \\
                                                                                                                                               1& 2^k & 2^{k-1}(2^k+1)-1 & 2^{k-1}(2^k-1)-1 & 2\\
                                                                                                                                            \end{array}
                                                                                                                                          \right)
\end{eqnarray*}

Thus, $\overline{\Gamma_{cent}(PSL(2,2^k))}$ is Laplacian integral for all values of $k.$ The following theorem describes the signless Laplacian spectrum of $\overline{\Gamma_{cent}(PSL(2,2^k))}.$

\bt\label{dhhhqt1} Let $\overline{\Gamma_{cent}(PSL(2,2^k))}$ be the co-centralizer graph of the projective special linear group $PSL(2,2^k)$. Then its signless Laplacian spectrum consists of:
\begin{enumerate}
\item[{\rm(a)}] $2^{2k}$ with multiplicity $2^k,$
\item[{\rm(b)}] $(2^{k-1}+2^{2k-1}+1)$ with multiplicity $2^{k-1}(2^k+1)-1,$
\item[{\rm(c)}] $3\times2^{k-1}+2^{2k-1}+1$  with multiplicity $2^{k-1}(2^k-1)-1,$ and
\item[{\rm (d)}] the three eigenvalues of the matrix $$\mathfrak{L}_P=\left[
                                                        \begin{array}{c|c|c}
                                                          2^{2k} & 2^{k-1}(2^k+1) & 2^{k-1}(2^k-1) \\
                                                          \hline
                                                          2^k+1 & 2^{k-1}+2^{2k-1}+1 & 2^{k-1}(2^k-1) \\
                                                          \hline
                                                          2^k+1 & 2^{k-1}(2^k+1) &  3\times 2^{k-1}+2^{2k-1}+1 \\
                                                        \end{array}
                                                      \right].
$$
\end{enumerate}\et
{\bf Proof.} With a suitable labeling of the vertices, the  signless Laplacian matrix for $\overline{\Gamma_{cent}(PSL(2,2^k))}$ can be written as
$$Q(\overline{\Gamma_{cent}(PSL(2,2^k))})=\left[
                                              \begin{array}{c|c|c}
                                                2^{2k}I & J & J \\
                                                \hline
                                                J & (2^{k-1}+2^{2k-1}+1)I & J \\
                                                \hline
                                                J & J & (3(2^{k-1})+2^{2k-1}+1)I \\
                                              \end{array}
                                            \right]
.$$

Now, $Q(\overline{\Gamma_{cent}(PSL(2,2^k))})\left(
           \begin{array}{c}
             -1\\
             1 \\
             \0_{2^k-1}\\
             \hline
             \0_{2^{k-1}(2^k+1)}\\
             \hline
             \0_{2^{k-1}(2^k-1)}\\
           \end{array}
         \right)=2^{2k}\left(
           \begin{array}{c}
             -1\\
             1 \\
             \0_{2^k-1}\\
             \hline
             \0_{2^{k-1}(2^k+1)}\\
             \hline
             \0_{2^{k-1}(2^k-1)}\\
           \end{array}
         \right).$

Therefore, $2^{2k}$ is an eigenvalue of $Q(\overline{\Gamma_{cent}(PSL(2,2^k))})$ with the following set $S_1$ of $2^k$ linearly independent eigenvectors;

$S_1= \left\{\left(
          \begin{array}{c}
            -1 \\
            1 \\
            0_{2^k-1} \\
            \hline
            \0_{2^{k-1}(2^k+1)} \\
            \hline
            \0_{2^{k-1}(2^k-1)} \\
          \end{array}
        \right),\left(
                  \begin{array}{c}
                    -1 \\
                    0\\
                    1 \\
                    \0_{2^k-2} \\
                    \hline
                    \0_{2^{k-1}(2^k+1)} \\
            \hline
            \0_{2^{k-1}(2^k-1)} \\
                  \end{array}
                \right), \ldots, \left(
                                   \begin{array}{c}
                                     -1 \\
                                     \0_{2^k-1} \\
                                     1 \\
                                     \hline
                                     \0_{2^{k-1}(2^k+1)} \\
            \hline
            \0_{2^{k-1}(2^k-1)} \\
                                   \end{array}
                                 \right)\right\}$.

Similarly, $$Q(\overline{\Gamma_{cent}(PSL(2,2^k))})\left(
                 \begin{array}{c}
                   \0_{2^k+1}\\
                   \hline
                   -1 \\
                   1 \\
                   \0_{2^{k-1}(2^k+1)-2} \\
                   \hline
                   \0_{2^{k-1}(2^k-1)} \\
                 \end{array}
               \right)=(2^{k-1}+2^{2k-1}+1)\left(
                               \begin{array}{c}
                                \0_{2^k+1}\\
                   \hline
                   -1 \\
                   1 \\
                   \0_{2^{k-1}(2^k+1)-2} \\
                   \hline
                   \0_{2^{k-1}(2^k-1)} \\
                               \end{array}
                             \right),$$
shows that  $(2^{k-1}+2^{2k-1}+1)$ is an eigenvalue of $Q(\overline{\Gamma_{cent}(PSL(2,2^k))})$ and in this way we can construct the following set $S_2$ of $2^{k-1}(2^k+1)-1$ linearly independent eigenvectors corresponding to $(2^{k-1}+2^{2k-1}+1);$
$$S_2= \left\{\left(
                 \begin{array}{c}
                   \0_{2^k+1}\\
                   \hline
                   -1 \\
                   1 \\
                   \0_{2^{k-1}(2^k+1)-2} \\
                   \hline
                   \0_{2^{k-1}(2^k-1)} \\
                 \end{array}
               \right),\left(
                 \begin{array}{c}
                   \0_{2^k+1}\\
                   \hline
                   -1 \\
                   0 \\
                   1 \\
                   \0_{2^{k-1}(2^k+1)-3} \\
                   \hline
                   \0_{2^{k-1}(2^k-1)} \\
                 \end{array}
               \right), \ldots,\left(
                 \begin{array}{c}
                   \0_{2^k+1}\\
                   \hline
                   -1 \\
                   \0_{2^{k-1}(2^k+1)-2} \\
                   1 \\
                   \hline
                   \0_{2^{k-1}(2^k-1)} \\
                 \end{array}
               \right)\right\}.$$

Moreover, $$Q(\overline{\Gamma_{cent}(PSL(2,2^k))})\left(
                 \begin{array}{c}
                   \0_{2^k+1}\\
                   \hline
                   \0_{2^{k-1}(2^k+1)} \\
                   \hline
                   -1\\
                   1\\
                   \0_{2^{k-1}(2^k-1)-2} \\
                 \end{array}
               \right)=(3(2^{k-1})+2^{2k-1}+1)\left(
                 \begin{array}{c}
                   \0_{2^k+1}\\
                   \hline
                   \0_{2^{k-1}(2^k+1)} \\
                   \hline
                   -1 \\
                   1 \\
                   \0_{2^{k-1}(2^k-1)-2} \\
                 \end{array}
               \right).$$
So, $(3(2^{k-1})+2^{2k-1}+1)$ is an eigenvalue of $Q(\overline{\Gamma_{cent}(PSL(2,2^k))}),$ and the following set $S_3$ lists $2^{k-1}(2^k-1)-1$ independent eigenvectors corresponding to this eigenvalue;
 $$S_3= \left\{\left(
                 \begin{array}{c}
                   \0_{2^k+1}\\
                   \hline
                   \0_{2^{k-1}(2^k+1)} \\
                   \hline
                   -1 \\
                   1 \\
                   \0_{2^{k-1}(2^k-1)-2} \\
                 \end{array}
               \right),\left(
                 \begin{array}{c}
                   \0_{2^k+1}\\
                   \hline
                   \0_{2^{k-1}(2^k+1)} \\
                   \hline
                   -1 \\
                   0 \\
                   1 \\
                   \0_{2^{k-1}(2^k-1)-3} \\
                 \end{array}
               \right), \ldots,\left(
                 \begin{array}{c}
                   \0_{2^k+1}\\
                   \hline
                   \0_{2^{k-1}(2^k+1)} \\
                   \hline
                   -1 \\
                   \0_{2^{k-1}(2^k-1)-2} \\
                   1 \\
                 \end{array}
               \right)\right\}.$$

Thus, we have obtained $2^k+2^{k-1}(2^k+1)-1+2^{k-1}(2^k-1)-1=2^k+2^{2k}-2$ eigenvalues of $Q(\overline{\Gamma_{cent}(PSL(2,2^k))}).$ Moreover, we note that all the eigenvectors constructed so far,
are orthogonal to $\left[
                     \begin{array}{c}
                       \1_{2^k+1} \\
                       \hline
                       \0_{2^{k-1}(2^k+1)} \\
                       \hline
                       \0_{2^{k-1}(2^k-1)}\\
                     \end{array}
                   \right],~\left[
                     \begin{array}{c}
                       \0_{2^k+1} \\
                       \hline
                       \1_{2^{k-1}(2^k+1)} \\
                       \hline
                       \0_{2^{k-1}(2^k-1)}\\
                     \end{array}
                   \right]$ and $\left[
                     \begin{array}{c}
                       \0_{2^k+1} \\
                       \hline
                       \0_{2^{k-1}(2^k+1)} \\
                       \hline
                       \1_{2^{k-1}(2^k-1)}\\
                     \end{array}
                   \right].$
Therefore, these three vectors span the remaining three eigenvectors of $Q(\overline{\Gamma_{cent}(PSL(2,2^k))}).$ Thus, the remaining eigenvectors
of $Q(\overline{\Gamma_{cent}(PSL(2,2^k))})$ are of the form $\left[
                     \begin{array}{c}
                       a\1_{2^k+1} \\
                       \hline
                       b\1_{2^{k-1}(2^k+1)} \\
                       \hline
                       c\1_{2^{k-1}(2^k-1)}\\
                     \end{array}
                   \right],$
for some $(a,b,c)\ne(0,0,0).$ Therefore if $\mu$ is an eigenvalue of $Q(\overline{\Gamma_{cent}(PSL(2,2^k))})$ with eigenvector $\left[
                     \begin{array}{c}
                       a\1_{2^k+1} \\
                       \hline
                       b\1_{2^{k-1}(2^k+1)} \\
                       \hline
                       c\1_{2^{k-1}(2^k-1)}\\
                     \end{array}
                   \right],$
then $a,b,c$ are the solution of the following system of equation
\begin{eqnarray*}
  \left(2^{2k}\right)a+\left(2^{k-1}\times (2^k+1)\right)b+\left(2^{k-1}\times(2^k-1)\right)c &=& 0 \\
   (2^k+1)a +\left(2^{k-1}+2^{2k-1}+1\right)b+\left(2^{k-1}\times(2^k-1)\right)c&=& 0 \\
    (2^k+1)a+\left(2^{k-1}\times(2^k+1)\right)b+\left(3\times2^{k-1}+2^{2k-1}+1\right)c&=& 0.
\end{eqnarray*}
Therefore, the remaining three eigenvalues of $Q(\overline{\Gamma_{cent}(PSL(2,2^k))})$ are the eigenvalues of the matrix
$\mathfrak{L}_P.$\qed

Hence, by Theorem~\ref{dhhhqt1}, $\overline{\Gamma_{cent}(PSL(2,2^k))}$ is signless Laplacian integral if $\mathfrak{L}_P$ have integral spectrum.
\section{Conclusion}
In this article, we have investigated the adjacency, (signless) Laplacian spectra of centralizer and co-centralizer graphs of the generalized quaternion group, the dihedral group, the quasidihedral group, the metacyclic group, and the projective special linear group. We also obtain conditions under which these graphs will be adjacency, (signless) Laplacian integral.

 \end{document}